\documentclass[12pt]{amsart}
\newtheorem{proposition}{Proposition}
\newtheorem{corollary}{Corollary}
\newtheorem{theorem}{Theorem}
\newtheorem{lemma}{Lemma}
\begin{document}
\title{}
\author{}
\maketitle

\centerline{\Large On the boundary behavior of the holomorphic}

\vskip 0.1in

\centerline{\Large sectional curvature of the Bergman metric}

\vskip 1cm

\centerline{\large Elisabetta Barletta\footnote{Universit\`a degli
Studi della Basilicata, Dipartimento di Matematica, Contrada
Macchia Romana, 85100 Potenza, Italy, e-mail: {\tt
barletta@unibas.it}}}

\begin{abstract} We obtain a conceptually new differential geometric proof of P.F.
Klembeck's result (cf. \cite{kn:Klem}) that the holomorphic
sectional curvature $k_g (z)$ of the Bergman metric of a strictly
pseudoconvex domain $\Omega \subset {\mathbb C}^n$ approaches
$-4/(n+1)$ (the constant sectional curvature of the Bergman metric
of the unit ball) as $z \to \partial \Omega$.
\end{abstract}

\section{Introduction}
Given a smoothly bounded strictly pseudoconvex domain $\Omega
\subset {\mathbb C}^n$ C.R. Graham \& J.M. Lee studied (cf.
\cite{kn:GrLe}) the $C^\infty$ regularity up to the boundary for
the solution to the Dirichlet problem $\Delta_g u = 0$ in $\Omega$
and $u = f$ on $\partial \Omega$, where $\Delta_g$ is the
Laplace-Beltrami operator of the Bergman metric $g$ of $\Omega$.
If $\varphi \in C^\infty (U)$ is a defining function ($\Omega = \{
z \in U : \varphi (z) < 0 \}$) their approach is to consider the
foliation $\mathcal F$ of a one-sided neighborhood $V$ of the
boundary $\partial \Omega$ by level sets $M_\epsilon = \{ z \in V
: \varphi (z) = - \epsilon \}$ ($\epsilon > 0$). Then $\mathcal F$
is a tangential CR foliation (cf. S. Dragomir \& S. Nishikawa,
\cite{kn:DrNi}) each of whose leaves is strictly pseudoconvex and
one may express $\Delta_g u = 0$ in terms of pseudohermitian
invariants of the leaves and the transverse curvature $r = 2 \,
\partial \overline{\partial} \varphi (\xi , \overline{\xi})$ and
its derivatives (the meaning of $\xi$ is explained in the next
section). The main technical ingredient is an ambient linear
connection $\nabla$ on $V$ whose pointwise restriction to each
leaf of $\mathcal F$ is the Tanaka-Webster connection (cf. S.
Webster, \cite{kn:Web}, and N. Tanaka, \cite{kn:Tan}) of the leaf.
An axiomatic description (and index free proof) of the existence
and uniqueness of $\nabla$ (referred to as the {\em Graham-Lee
connection} of $(V , \varphi )$) was provided in \cite{kn:BDU}. As
a natural continuation of the ideas in \cite{kn:GrLe} one may
relate the Levi-Civita connection $\nabla^g$ of $(V , g)$ to the
Graham-Lee connection $\nabla$ and compute the curvature $R^g$ of
$\nabla^g$ in terms of the curvature of $\nabla$. Together with an
elementary asymptotic analysis (as $\epsilon \to 0$) this leads to
a purely differential geometric proof of the result of P.F.
Klembeck, \cite{kn:Klem}, that the sectional curvature of $(\Omega
, g)$ tends to $-4/(n+1)$ near the boundary $\partial \Omega$. The
Author believes that one cannot overestimate the importance of the
Graham-Lee connection (and that the identities (\ref{e:425}) and
(\ref{e:434}) in Section 3 admit other applications as well, e.g.
in the study of the geometry of the second fundamental form of a
submanifold in $(\Omega , g)$).

\section{The Levi-Civita versus the Graham-Lee connection}
Let $\Omega$ be a smoothly bounded strictly pseudoconvex domain in
${\mathbb C}^n$ and $K(z, \zeta )$ its Bergman kernel (cf. e.g.
\cite{kn:Hel}, p. 364-371). As a simple application of C.
Fefferman's asymptotic development (cf. \cite{kn:Fef}) of the
Bergman kernel  $\varphi (z) = - K(z,z)^{-1/(n+1)}$ is a defining
function for $\Omega$ (and $\Omega = \{ \varphi < 0 \}$). Cf. A.
Kor\'anyi \& H.M. Reimann, \cite{kn:KoRe}, for a proof. Let us set
$\theta = \frac{i}{2}(\overline{\partial} -
\partial )\varphi$.
Then $d \theta = i \,
\partial \overline{\partial} \varphi$. Let us differentiate $\log
|\varphi | = - (1/(n+1)) \log K$ (where $K$ is short for $K(z,z)$)
so that to obtain
\[ \frac{1}{\varphi} \, \overline{\partial} \varphi = -
\frac{1}{n+1} \, \overline{\partial} \log K. \] Applying the
operator $i \, \partial$ leads to
\begin{equation} \frac{1}{\varphi} \; d \theta -
\frac{i}{\varphi^2} \; \partial \varphi \wedge \overline{\partial}
\varphi = - \frac{i}{n+1} \, \partial \overline{\partial} \log K.
\label{e:b3}
\end{equation}
We shall need the Bergman metric $g_{j\overline{k}} = \partial^2
\log K /\partial z^j \partial \overline{z}^k$. This is well known
to be a K\"ahler metric on $\Omega$.
\begin{proposition} For any smoothly bounded strictly pseudoconvex
domain $\Omega \subset {\mathbb C}^n$ the Bergman metric $g$ is
given by
\begin{equation}
g(X,Y) = \frac{n+1}{\varphi} \{ \frac{i}{\varphi} \, (\partial
\varphi \wedge \overline{\partial} \varphi )(X , J Y ) - d \theta
(X , J Y) \} , \label{e:b4}
\end{equation}
for any $X,Y \in \mathcal{X}(\Omega )$.
\end{proposition}
{\em Proof}. Let $\omega (X,Y) = g(X , J Y)$ be the K\"ahler
$2$-form of $(\Omega ,J, g)$, where $J$ is the underlying complex
structure. Then $\omega = - i \,
\partial \overline{\partial} \log K$ and (\ref{e:b3}) may be
written in the form (\ref{e:b4}). Q.e.d.

\vskip 0.1in We denote by $M_\epsilon = \{ z \in \Omega : \varphi
(z) = - \epsilon \}$ the level sets of $\varphi$. For $\epsilon >
0$ sufficiently small $M_\epsilon$ is a strictly pseudoconvex CR
manifold (of CR dimension $n-1$). Therefore, there is a one-sided
neighborhood $V$ of $\partial \Omega$ which is foliated by the
level sets of $\varphi$. Let $\mathcal{F}$ be the relevant
foliation and let us denote by $H(\mathcal{F}) \to V$
(respectively by $T_{1,0}(\mathcal{F}) \to V$) the bundle whose
portion over $M_\epsilon$ is the Levi distribution $H(M_\epsilon
)$ (respectively the CR structure $T_{1,0}(M_\epsilon )$) of
$M_\epsilon$. Note that
\[ T_{1,0}(\mathcal{F}) \cap T_{0,1}(\mathcal{F}) = (0), \]
\[ [\Gamma^\infty (T_{1,0}(\mathcal{F})) , \Gamma^\infty (T_{1,0}(\mathcal{F}))]
\subseteq \Gamma^\infty (T_{1,0}(\mathcal{F})). \] Here
$T_{0,1}(\mathcal{F}) = \overline{T_{1,0}(\mathcal{F})}$. For a
review of the basic notions of CR and pseudohermitian geometry
needed through this paper one may see S. Dragomir \& G. Tomassini,
\cite{kn:DrTo}. Cf. also S. Dragomir, \cite{kn:Dra}. By a result
of J.M. Lee \& R. Melrose, \cite{kn:LeMe}, there is a unique
complex vector field $\xi$ on $V$, of type $(1,0)$, such that
$\partial \varphi (\xi ) = 1$ and $\xi$ is orthogonal to
$T_{1,0}(\mathcal{F})$ with respect to $\partial
\overline{\partial} \varphi$ i.e. $\partial \overline{\partial}
\varphi (\xi , \overline{Z}) = 0$ for any $Z \in
T_{1,0}(\mathcal{F})$. Let $r = 2 \, \partial \overline{\partial}
\varphi (\xi , \overline{\xi})$ be the {\em transverse curvature}
of $\varphi$. Moreover let $\xi = \frac{1}{2}(N - i T)$ be the
real and imaginary parts of $\xi$. Then
\[ (d \varphi )(N) = 2, \;\;\; (d \varphi )(T) = 0, \]
\[ \theta (N) = 0, \;\;\; \theta (T) = 1, \]
\[ \partial \varphi (N) = 1, \;\;\; \partial \varphi (T) = i. \]
In particular $T$ is tangent to (the leaves of) $\mathcal{F}$. Let
$g_\theta$ be the tensor field given by
\begin{equation}
g_\theta (X,Y) = (d \theta )(X, J Y), \;\; g_\theta (X,T) = 0,
\;\; g_\theta (T , T) = 1, \label{e:A.1}
\end{equation}
for any $X,Y \in H(\mathcal{F})$. Then $g_\theta$ is a tangential
Riemannian metric for $\mathcal{F}$ i.e. a Riemannian metric in
$T(\mathcal{F}) \to V$. Note that the pullback of $g_\theta$ to
each leaf $M_\epsilon$ of $\mathcal{F}$ is the Webster metric of
$M_\epsilon$ (associated to the contact form $j_\epsilon^*
\theta$, where $j_\epsilon : M_\epsilon \subset V$). As a
consequence of (\ref{e:b4}), $J T = -N$ and $i_N \, d \theta = r
\, \theta$ (see also (\ref{e:A.4}) below)
\begin{corollary} The Bergman metric $g$ of $\Omega \subset
{\mathbb C}^n$ is given by
\begin{equation}
g(X,Y) = - \frac{n+1}{\varphi} \, g_\theta (X,Y), \;\;\; X,Y \in
H(\mathcal{F}). \label{e:b5}
\end{equation}
\begin{equation}
g(X,T) = 0, \;\; g(X, N) = 0, \;\;\; X \in H(\mathcal{F}),
\label{e:b6}
\end{equation}
\begin{equation}
g(T , N ) = 0, \;\; g(T,T) = g(N,N) = \frac{n+1}{\varphi} \left(
\frac{1}{\varphi} - r \right) . \label{e:b7}
\end{equation}
In particular $1 - r \varphi > 0$ everywhere in $\Omega$.
\label{c:1}
\end{corollary}
Using (\ref{e:b5})-(\ref{e:b7}) we may relate the Levi-Civita
connection $\nabla^g$ of $(V , g)$ to another canonical linear
connection on $V$, namely the {\em Graham-Lee connection} of
$\Omega$. The latter has the advantage of staying finite at the
boundary (it gives the Tanaka-Webster connection of $\partial
\Omega$ as $z \to \partial \Omega$). We proceed to recalling the
Graham-Lee connection. Let $\{ W_\alpha : 1 \leq \alpha \leq n-1
\}$ be a local frame of $T_{1,0}(\mathcal{F})$, so that $\{
W_\alpha , \xi \}$ is a local frame of $T^{1,0}(V)$. We consider
as well
\[ L_\theta (Z,  \overline{W}) \equiv - i (d \theta )(Z, \overline{W}), \;\;\; Z,W
\in T_{1,0}(\mathcal{F}). \] Note that $L_\theta$ and (the
$\mathbb{C}$-linear extension of) $g_\theta$ coincide on
$T_{1,0}(\mathcal{F}) \otimes T_{0,1}(\mathcal{F})$. We set
$g_{\alpha\overline{\beta}} = g_\theta (W_\alpha ,
W_{\overline{\beta}})$. Let $\{ \theta^\alpha : 1 \leq \alpha \leq
n-1 \}$ be the (locally defined) complex $1$-forms on $V$
determined by \[ \theta^\alpha (W_\beta ) = \delta^\alpha_\beta \,
, \;\; \theta^\alpha (W_{\overline{\beta}}) = 0, \;\;
\theta^\alpha (T) = 0, \;\; \theta^\alpha (N) = 0. \] Then $\{
\theta^\alpha , \, \theta^{\overline{\alpha}} , \, \theta , \, d
\varphi \}$ is a local frame of $T(V) \otimes \mathbb{C}$ and one
may easily show that
\begin{equation}
d \theta = 2 i g_{\alpha\overline{\beta}} \, \theta^\alpha \wedge
\theta^{\overline{\beta}} + r \, d \varphi \wedge \theta .
\label{e:A.2}
\end{equation}
As an immediate consequence
\begin{equation}
i_T \, d \theta = - \frac{r}{2} \, d \varphi , \;\;\; i_N \, d
\theta = r \, \theta . \label{e:A.4}
\end{equation}
As an application of (\ref{e:A.2}) we decompose $[T,N]$ (according
to $T(V) \otimes \mathbb{C} = T_{1,0}(\mathcal{F}) \oplus
T_{0,1}(\mathcal{F}) \oplus \mathbb{C} T \oplus \mathbb{C} N$) and
obtain
\begin{equation}
[T,N] = i \, W^\alpha (r) W_\alpha - i \, W^{\overline{\alpha}}(r)
W_{\overline{\alpha}} + 2 r T, \label{e:A.5}
\end{equation}
where $W^\alpha (r) = g^{\alpha\overline{\beta}}
W_{\overline{\beta}}(r)$ and $W^{\overline{\alpha}} (r) =
\overline{W^\alpha (r)}$.
\par
Let $\nabla$ be a linear connection on $V$. Let us consider the
$T(V)$-valued $1$-form $\tau$ on $V$ defined by
\[ \tau (X) = T_\nabla (T , X), \;\;\; X \in T(V), \]
where $T_\nabla$ is the torsion tensor field of $\nabla$. We say
$T_\nabla$ is {\em pure} if
\begin{equation}
T_\nabla (Z, W) = 0, \;\; T_\nabla (Z, \overline{W}) = 2 i
L_\theta (Z , \overline{W}) T, \label{e:A.6}
\end{equation}
\begin{equation}
T_\nabla (N , W) = r \, W + i \, \tau (W), \label{e:A.7}
\end{equation}
for any $Z,W \in T_{1,0}(\mathcal{F})$, and
\begin{equation}
\tau (T_{1,0}(\mathcal{F})) \subseteq T_{0,1}(\mathcal{F}),
\label{e:A.8}
\end{equation}
\begin{equation}
\tau (N) = -  \, J \, \nabla^H r -  2 r \, T. \label{e:A.9}
\end{equation}
Here $\nabla^H r$ is defined by $\nabla^H r = \pi_H \nabla r$ and
$g_\theta (\nabla r , X) = X(r)$, $X \in T(\mathcal{F})$. Also
$\pi_H : T(\mathcal{F}) \to H(\mathcal{F})$ is the projection
associated to the direct sum decomposition $T(\mathcal{F}) =
H(\mathcal{F}) \oplus \mathbb{R} T$. We recall the following
\begin{theorem} There is a unique linear connection $\nabla$ on $V$ such that {\rm
i)} $T_{1,0}(\mathcal{F})$ is parallel with respect to $\nabla$,
{\rm ii)} $\nabla L_\theta = 0$, $\nabla T = 0$, $\nabla N = 0$,
and {\rm iii)} $T_\nabla$ is pure. \label{t:A.1}
\end{theorem} \noindent  $\nabla$ given by
Theorem \ref{t:A.1} is the {\em Graham-Lee connection}. Theorem
\ref{t:A.1} is essentially Proposition 1.1 in \cite{kn:GrLe}, p.
701-702. The axiomatic description in Theorem \ref{t:A.1} is due
to \cite{kn:DrNi} (cf. Theorem 2 there). An index-free proof of
Theorem \ref{t:A.1} was given in \cite{kn:BDU} relying on the
following
\begin{lemma} Let $\phi : T(\mathcal{F}) \to T(\mathcal{F})$ be
the bundle morphism given by $\phi (X) = J X$, for any $X \in
H(\mathcal{F})$, and $\phi (T) = 0$. Then
\[ \phi^2 = - I + \theta \otimes T, \]
\[ g_\theta (X , T) = \theta (X), \]
\[ g_\theta (\phi X , \phi Y) = g_\theta (X,Y) - \theta (X) \theta (Y), \]
for any $X,Y \in T(\mathcal{F})$. Moreover, if $\nabla$ is a
linear connection on $V$ satisfying the axioms {\rm (i)-(iii)} in
Theorem $\ref{t:A.1}$ then
\begin{equation}
\phi \circ \tau + \tau \circ \phi = 0 \label{e:A.10}
\end{equation}
along $T(\mathcal{F})$. Consequently $\tau$ may be computed as
\begin{equation}
\tau (X) = - \frac{1}{2} \phi (\mathcal{L}_T \phi ) X,
\label{e:A.11}
\end{equation} for any $X \in H(\mathcal{F})$.
\label{l:A.1}
\end{lemma}
A rather lengthy but straightforward calculation (based on
Corollary \ref{c:1}) leads to
\begin{theorem} Let $\Omega \subset {\mathbb C}^n$ be a smoothly bounded strictly
pseudoconvex domain, $K(z, \zeta )$ its Bergman kernel, and
$\varphi (z) = - K(z,z)^{-1/(n+1)}$. Then the Levi-Civita
connection $\nabla^g$ of the Bergman metric and the Graham-Lee
connection of $(\Omega , \varphi )$ are related by
\begin{equation} \nabla^g_X Y = \nabla_X Y +
\label{e:b13}
\end{equation}
\[ + \left\{ \frac{\varphi}{1  - \varphi r} \, g_\theta (\tau X ,
Y) + g_\theta (X , \phi Y) \right\} T - \] \[ - \left\{ g_\theta
(X,Y) + \frac{\varphi}{1  - \varphi r} \, g_\theta (X , \phi \,
\tau \, Y) \right\} N, \]
\begin{equation}
\nabla^g_X T = \tau X - \left( \frac{1}{\varphi} - r \right) \phi
X -  \label{e:b17}
\end{equation}
\[ - \frac{\varphi}{2(1 - r \varphi )} \left\{ X(r) T + (\phi X)(r)
N \right\} , \]
\begin{equation}
\nabla^g_X N = - \left( \frac{1}{\varphi} - r \right) X + \tau \,
\phi\, X +  \label{e:b21}
\end{equation}
\[ + \frac{\varphi}{2(1 - r \varphi )} \{ (\phi X)(r) \, T -
X(r) \, N \} , \]
\begin{equation}
\nabla^g_T X = \nabla_T X - \left( \frac{1}{\varphi} - r \right)
\phi X -  \label{e:b25}
\end{equation}
\[ - \frac{\varphi}{2(1 - r \varphi )} \{ X(r) T + (\phi X)(r)
N \} , \]
\begin{equation}
\nabla^g_N X = \nabla_N X - \frac{1}{\varphi} \, X +
\label{e:b29}
\end{equation}
\[ + \frac{\varphi}{2(1 - r \varphi )} \{ (\phi X)(r) T - X(r) N \}
, \]
\begin{equation}
\nabla^g_N T = - \frac{1}{2} \, \phi \, \nabla^H r -
 \label{e:b30}
\end{equation}
\[ - \frac{\varphi}{2(1 - r \varphi )} \left\{ \left( N(r) +
\frac{4}{\varphi^2} - \frac{2 r}{\varphi} \right) T + T(r) N
\right\} . \]
\begin{equation}
\nabla^g_T N = \frac{1}{2} \, \phi \nabla^H r - \label{e:b31}
\end{equation}
\[ - \frac{\varphi}{2(1 - r \varphi )} \left\{ \left( N(r) +
\frac{4}{\varphi^2} - \frac{6 r}{\varphi} + 4 r^2 \right) T + T(r)
N \right\} , \]
\begin{equation}
\nabla^g_T T = - \frac{1}{2} \; \nabla^H r - \label{e:b32}
\end{equation}
\[ - \frac{\varphi}{2(1 - r \varphi )} \left\{ T(r) T - \left(
N(r) + \frac{4}{\varphi^2} - \frac{6 r}{\varphi} + 4 r^2 \right) N
\right\} , \]
\begin{equation}
\nabla^g_N N = - \frac{1}{2} \; \nabla^H r + \label{e:b33}
\end{equation}
\[ + \frac{\varphi}{2(1 - r \varphi )} \left\{ T(r) T - \left(
N(r) + \frac{4}{\varphi^2} - \frac{2 r}{\varphi} \right) N
\right\} , \] for any $X, Y \in H({\mathcal F})$.
\end{theorem}

\section{Klembeck's theorem}
The original proof of the result by P.F. Klembeck (cf. Theorem 1
in \cite{kn:Klem}, p. 276) employs a formula of S. Kobayashi,
\cite{kn:Kob}, expressing the components
$R_{j\overline{k}r\overline{s}}$ of the Riemann-Christoffel
$4$-tensor of $(\Omega , g)$ as
\[ - \frac{1}{2} R_{j\overline{k}r\overline{s}} =
g_{j\overline{k}} g_{r\overline{s}} + g_{j\overline{s}}
g_{r\overline{k}} - \frac{1}{K^2} \{ K \,
K_{j\overline{k}r\overline{s}} - K_{jr} K_{\overline{k} \,
\overline{s}} \} + \]
\[ + \frac{1}{K^4} \sum_{\ell , m} g^{\overline{\ell}m} \{ K \,
K_{jr\overline{\ell}} - K_{jr} K_{\overline{\ell}} \} \{ K \,
K_{\overline{k} \, \overline{s} m} - K_{\overline{k} \,
\overline{s}} K_m \} \] where $K = K(z,z)$ and its indices denote
derivatives. However the calculation of the inverse matrix
$[g^{j\overline{k}}] = [g_{j\overline{k}}]^{-1}$ turns out to be a
difficult problem and \cite{kn:Klem} only provides an asymptotic
formula as $z \to \partial \Omega$. Our approach is to compute the
holomorphic sectional curvature of $(\Omega , g)$ by deriving an
explicit relation among the curvature tensor fields $R^g$ and $R$
of the Levi-Civita and Graham-Lee connections respectively. We
start by recalling a pseudohermitian analog to holomorphic
curvature (built by S.M. Webster, \cite{kn:Web}).
\par
Let $M$ be a nondegenerate CR manifold of type $(n-1,1)$ and
$\theta$ a contact form on $M$. Let $G_1 (H(M))_x$ consist of all
$2$-planes $\sigma \subset T_x (M)$ such that i) $\sigma \subset
H(M)_x$ and ii) $J_x (\sigma ) = \sigma$. Then $G_1 (H(M))$ (the
disjoint union of all $G_1 (H(M))_x$) is a fibre bundle over $M$
with standard fibre ${\mathbb C}P^{n-2}$. Let $R^\nabla$ be the
curvature of the Tanaka-Webster connection $\nabla$ of $(M ,
\theta )$. We define a function $k_\theta : G_1 (H(M)) \rightarrow
{\mathbb R}$ by setting
\[
k_\theta (\sigma ) = - \frac{1}{4} R^\nabla_x (X , J_x X , X , J_x
X)
\] for any $\sigma \in G_1 (H(M))$ and any linear basis $\{ X ,
J_x X \}$ in $\sigma$ satisfying $G_\theta (X , X) = 1$. It is a
simple matter that the definition of $k_\theta (\sigma )$ does not
depend upon the choice of orthonormal basis $\{ X , J_x X \}$, as
a consequence of the following properties
\[ R^\nabla (Z,W,X,Y) + R^\nabla (Z,W,Y,X) = 0, \]
\[ R^\nabla (Z,W,X,Y) + R^\nabla (W,Z,X,Y) = 0. \]
 $k_\theta$ is referred to as the
({\em pseudohermitian}) {\em sectional curvature} of $(M , \theta
)$. As mentioned above the notion is due to S.M. Webster,
\cite{kn:Web}, who also gave examples of pseudohermitian space
forms (pseudohermitian manifolds $(M , \theta )$ with $k_\theta$
constant). Cf. also \cite{kn:BaDr2} for a further study of contact
forms of constant pseudohermitian sectional curvature. With
respect to an arbitrary (not necessarily orthonormal) basis $\{ X
, J_x X \}$ of the $2$-plane $\sigma$ the sectional curvature
$k_\theta (\sigma )$ is also expressed by
\[ k_\theta (\sigma ) = - \frac{1}{4}
\frac{R^\nabla_x (X , J_x X , X , J_x X )}{G_\theta (X,X)^2}\, .
\] To prove this statement one merely applies the definition of
$k_\theta (\sigma )$ for the orthonormal basis $\{ U , J_x U \}$,
with $U = G_\theta (X,X)^{-1/2} X$.  As $X \in H(M)_x$ there is $Z
\in T_{1,0}(M)_x$ such that $X = Z + \overline{Z}$. Thus
\[ k_\theta (\sigma ) = \frac{1}{4}
\frac{R_x (Z , \overline{Z} , Z , \overline{Z})}{g_\theta (Z ,
\overline{Z})^2} \, . \] The coefficient $1/4$ is chosen such that
the sphere $S^{2n-1} \subset {\mathbb C}^{n}$ has constant
curvature $+1$. Cf. \cite{kn:DrTo}, Chapter 1. With the notations
in Section 2 let us set $f = \varphi /(1 - \varphi r)$. Then
\[ X(f) = f^2 \; X(r), \;\;\; X \in T({\mathcal F}). \]
Let $R^g$ and $R$ be respectively the curvature tensor fields of
the linear connections $\nabla^g$ and $\nabla$ (the Graham-Lee
connection). For any $X,Y,Z \in H({\mathcal F})$ (by
(\ref{e:b13}))
\[ \nabla^g_X \nabla^g_Y Z = \nabla^g_X \left( \nabla_Y Z + \left\{
f \, g_\theta (\tau (Y), Z) + g_\theta (Y , \phi Z) \right\} T -
\right.  \] \[ \left. - \left\{ g_\theta (Y,Z) + f \, g_\theta (Y
, \phi \tau (Z)) \right\} N \right) = \] by $\nabla_Y Z \in
H({\mathcal F})$ together with (\ref{e:b13})
\[ = \nabla_X \nabla_Y Z + \left\{ f \, g_\theta (\tau (X) ,
\nabla_Y Z ) + g_\theta (X , \phi \nabla_Y Z) \right\} T - \]
\[ - \left\{ g_\theta (X , \nabla_Y Z ) + f \, g_\theta (X , \phi
\tau (\nabla_Y Z)) \right\} N + \]
\[ + \left\{ f \, g_\theta (\tau (Y) , Z) + g_\theta (Y , \phi Z)
\right\} \nabla^g_X T + \] \[ + \left\{ X(f) g_\theta (\tau (Y),
Z) + f \, X(g_\theta (\tau (Y), Z)) + X(g_\theta (Y, \phi Z))
\right\} T - \] \[ - \left\{ g_\theta (Y,Z) + f\, g_\theta (Y ,
\phi \tau (Z)) \right\} \nabla^g_X N +
\]
\[ - \left\{ X (g_\theta (Y,Z)) + X(f) g_\theta (Y, \phi \tau (Z))
+ f \, X(g_\theta (Y , \phi \tau (Z))) \right\} N = \] by
(\ref{e:b17}), (\ref{e:b21})
\[ = \nabla_X \nabla_Y Z + \left\{ X(\Omega (Y,Z)) + \Omega (X , \nabla_Y Z)
+ \right. \] \[ \left. + X(f) A(Y,Z) + f \left[ X(A(Y,Z)) + A(X\,
\nabla_Y Z)\right] \right\} T - \]
\[ - \left\{ X(g_\theta (Y,Z)) + g_\theta (X , \nabla_Y Z) +
\right. \]
\[ + \left. X(f) \Omega (Y, \tau (Z)) + f \left[ X(\Omega (Y, \tau
(Z))) + \Omega (X, \tau (\nabla_Y Z))\right] \right\} N + \]
\[ + \left\{ f \, A(Y,Z) + \Omega (Y,Z) \right\} \times \] \[ \times \left\{ \tau (X)
- \frac{1}{f} \, \phi X - \frac{f}{2} \left( X(r) T + (\phi X)(r)
N \right) \right\} - \]
\[ - \left\{ g_\theta (Y,Z) + f \, \Omega (Y , \tau (Z)) \right\}
\times \] \[ \times  \left\{ - \frac{1}{f} \, X + \tau (\phi X) +
\frac{f}{2} \left( (\phi X)(r) T - X(r) N \right) \right\} \]
where we have set as usual $A(X,Y) = g_\theta (\tau (X) , Y)$ and
$\Omega (X,Y) = g_\theta (X , \phi Y)$. We may conclude that
\begin{equation}
\nabla^g_X \nabla^g_Y Z = \nabla_X \nabla_Y Z + [f \, A(Y,Z) +
\Omega (Y,Z)] \left( \tau (X) - \frac{1}{f} \, \phi X \right)  +
\label{e:423}
\end{equation}
\[ + [g_\theta (Y,Z) + f \, \Omega (Y, \tau (Z))] \left( \frac{1}{f} \,
X - \tau (\phi X)\right) + \]
\[ + \left\{ X(\Omega (Y,Z)) + \Omega (X , \nabla_Y Z) + f \left[
X(A(Y,Z)) + A(X , \nabla_Y Z)\right] + \right. \]
\[ +  \frac{f}{2} \left[ X(r) (f \, A(Y,Z) - \Omega (Y,Z)) - \right. \]
\[ \left. \left. - (\phi X)(r) (g_\theta (Y,Z) + f \, \Omega (Y, \tau (Z))) \right]
\right\} T - \]
\[ - \left\{ X(g_\theta (Y,Z)) + g_\theta (X , \nabla_Y Z) + f \left[ X(\Omega (Y, \tau
(Z))) + \Omega (X , \tau(\nabla_Y Z)) \right] - \right. \]
\[  - \frac{f}{2} \left[ X(r) (g_\theta (Y,Z) - f \, \Omega
(Y, \tau (Z))) - \right. \] \[ \left. \left. - (\phi X)(r) (f \,
A(Y,Z) + \Omega (Y,Z)) \right] \right\} N  \] for any $X,Y,Z \in
H({\mathcal F})$. Next we use the decomposition $[X,Y] = \pi_H
[X,Y] + \theta ([X,Y]) T$ and (\ref{e:b13}), (\ref{e:b25}) to
calculate
\[ \nabla^g_{[X,Y]} Z = \nabla^g_{\pi_H [X,Y]} Z + \theta ([X,Y])
\nabla^g_T Z = \]
\[ = \nabla_{\pi_H [X,Y]} Z + \left\{ f \, g_\theta (\tau (\pi_H
[X,Y]), Z) + g_\theta (\pi_H [X,Y] , \phi Z) \right\} T - \]
\[ - \left\{ g_\theta (\pi_H [X,Y] , Z) + f \, g_\theta (\pi_H
[X,Y] , \phi \tau (Z)) \right\} N + \]
\[ + \theta ([X,Y]) \left\{ \nabla_T Z - \frac{1}{f} \, \phi Z -
\frac{f}{2} (Z(r) T + (\phi Z)(r) N ) \right\} \] so that (by
$\tau (T) = 0$)
\begin{equation} \nabla^g_{[X,Y]} Z = \nabla_{[X,Y]} Z -
\frac{1}{f} \, \theta ([X,Y]) \phi Z + \label{e:424}
\end{equation}
\[ + \left\{ f \, A([X,Y],Z) + \Omega ([X,Y], Z) - \frac{f}{2}
\theta ([X,Y]) Z(r) \right\} T - \]
\[ - \left\{ g_\theta ([X,Y], Z) + f \, \Omega ([X,Y], \tau (Z)) +
\frac{f}{2} \theta ([X,Y]) (\phi Z)(r) \right\} N \] for any
$X,Y,Z \in H({\mathcal F})$. Consequently by
(\ref{e:423})-(\ref{e:424}) (and by $\nabla g_\theta = 0$, $\nabla
\Omega = 0$) we may compute
\[ R^g (X,Y) Z = \nabla^g_X \nabla^g_Y Z - \nabla^g_Y \nabla^g_X Z
- \nabla^g_{[X,Y]} Z \] so that to obtain
\begin{equation}
R^g (X,Y) Z = R(X,Y)Z + \frac{1}{f} \, \theta ([X,Y]) \phi Z +
\label{e:425}
\end{equation}
\[ + (f \, A(Y,Z) + \Omega (Y,Z)) \left( \tau (X) - \frac{1}{f} \,
\phi X \right) - \] \[ - (f \, A(X,Z) + \Omega (X,Z))\left( \tau
(Y) - \frac{1}{f} \, \phi Y \right) + \]
\[ + (g_\theta (Y,Z) + f \, \Omega (Y,  \tau (Z)) \left( \frac{1}{f} \,
X - \tau (\phi X)) \right) - \]
\[ - (g_\theta (X,Z) + f \, \Omega (X , \tau (Z))) \left(
\frac{1}{f} \, Y - \tau (\phi Y) \right) + \]
\[ + \left\{ f \left[ (\nabla_X A)(Y,Z) - (\nabla_Y A)(X,Z)
\right] + \right. \]
\[ +  \frac{f}{2} [ X(r) (f \, A(Y,Z) - \Omega (Y,Z)) - Y(r) (f \,
A(X,Z) - \Omega (X,Z)) - \]
\[ - (\phi X)(r) (g_\theta (Y,Z) + f \, \Omega (Y , \tau (Z))) +
(\phi Y)(r)(g_\theta (X,Z) + f \, \Omega (X , \tau (Z))) + \] \[
\left. + Z(r) \theta ([X,Y]) ] \right\} T - \]
\[ - \left\{ f \left[ \Omega (Y, (\nabla_X \tau )Z) - \Omega (X ,
(\nabla_Y \tau )Z ) \right] - \right. \]
\[ - \frac{f}{2} [X(r)(g_\theta (Y,Z) - f \, \Omega (Y , \tau
(Z))) - Y(r) (g_\theta (X,Z) - f \, \Omega (X , \tau (Z))) - \]
\[ - (\phi X)(r) (f \, A(Y,Z) + \Omega (Y,Z)) + (\phi Y)(r) (f \,
A(X,Z) + \Omega (X,Z)) + \]
\[ \left. + (\phi Z)(r) \theta ( [X,Y]) ] \right\} N \]
for any $X,Y,Z \in H({\mathcal F})$. Let us take the inner product
of (\ref{e:425}) with $W \in H({\mathcal F})$ and use
(\ref{e:b5})-(\ref{e:b6}). We obtain
\[ g(R^g (X,Y)Z , W) = -  \frac{n+1}{\varphi} \{ g_\theta (R(X,Y)Z
, W) - \frac{1}{f} \, \theta ( [X,Y]) \Omega (Z,W) + \]
\[ + [f\, A(Y,Z) + \Omega (Y,Z)] [A(X,W) + \frac{1}{f} \, \Omega
(X,W)] - \]
\[ - [f \, A(X,Z) + \Omega (X,Z)] [A(Y,W) + \frac{1}{f} \, \Omega
(Y,W)] + \]
\[ + [g_\theta (Y,Z) + f \, \Omega (Y, \tau (Z))] [\frac{1}{f} \,
g_\theta (X,W) + \Omega (X , \tau (W))] - \]
\[ - [g_\theta (X,Z) + f \, \Omega (X , \tau (Z)) ] [\frac{1}{f}
\, g_\theta (Y,W) + \Omega (Y , \tau (W))] \} . \] In particular
for $Z = Y$ and $W = X$ (as $\Omega = - d \theta$)
\[ g(R^g (X,Y)Y , X) = - \frac{n+1}{\varphi} \{ g_\theta
(R(X,Y)Y,X) + \] \[ +  \frac{2}{f} \, \Omega (X,Y)^2 + f \, A(X,X)
A(Y,Y) -
\]
\[ - \frac{1}{f} [f^2 \, A(X,Y)^2 - \Omega (X,Y)^2 ] + \]
\[ + \frac{1}{f} [g_\theta (X,X) + f \, \Omega (X , \tau
(X))][g_\theta (Y,Y) + f \, \Omega (Y, \tau (Y))] - \]
\[ - \frac{1}{f} [g_\theta (X,Y) + f \, \Omega (X , \tau (Y))]^2
\} . \] Note that
\[ A(\phi X , \phi X) = g_\theta (\tau (\phi X) , \phi X) = -
g_\theta (\phi \tau X , \phi X) = - A(X,X), \]
\[ \Omega (\phi X ,
\tau (\phi X)) = g_\theta (\phi X , \phi \tau (\phi X)) = g_\theta
(X , \tau (\phi X)) = \] \[ = - g_\theta (X , \phi \tau (X)) = -
\Omega (X, \tau (X)), \] \[ \Omega (X , \tau (\phi X)) = g_\theta
(X , \phi \tau (\phi X)) = - g_\theta (X , \tau (\phi^2 X)) = \]
\[ =  g_\theta (X , \tau (X)) = A(X,X). \] Hence
\begin{equation} g(R^g (X, \phi X)\phi X , X) = - \frac{n+1}{\varphi} \{
g_\theta (R(X, \phi X)\phi X , X) + \label{e:426}
\end{equation} \[ +
\frac{4}{f} \, g_\theta (X,X)^2 - 2 f [A(X,X)^2 + A(X , \phi X)^2
] \} . \] Let $\sigma \subset T({\mathcal F})_z$ be the $2$-plane
spanned by $\{ X , \phi_z X \}$ for $X \in H({\mathcal F})_z$, $X
\neq 0$. By (\ref{e:b5}) if $Y = \phi_z X$ then
\[ g_z (X,X) g_z (Y,Y) - g_z (X,Y)^2 = \] \[ = \left( \frac{n+1}{\varphi
(z)} \right)^2 \{ g_{\theta , z} (X,X) g_{\theta , z}(Y,Y) -
g_{\theta , z}(X,Y) \} = \] \[ =  \left( \frac{n+1}{\varphi (z)}
\right)^2 g_{\theta , z} (X,X)^2 \] so that (by (\ref{e:426})) the
sectional curvature $k_g (\sigma )$ of the $2$-plane $\sigma$ is
expressed by (for $Y = \phi_z X$)
\[ k_g (\sigma ) = \frac{g_z (R^g_z (X , Y)Y , X)}{g_z (X,X)g_z (Y,Y) - g_z (X,Y)^2} = \]
\[ = - \frac{\varphi (z)}{n+1} \{ - 4 k_\theta (\sigma ) +
\frac{4}{f(z)} - 2 f(z) \frac{A_z (X,X)^2 + A_z (X , \phi_z
X)^2}{g_{\theta , z}(X,X)^2} \}  \] where $k_\theta$ restricted to
a leaf of $\mathcal F$ is the pseudohermitian sectional curvature
of the leaf. Note that $k_\theta$ and $A$ stay finite at the
boundary (and give respectively the pseudohermitian sectional
curvature and the pseudohermitian torsion of $(\partial \Omega ,
\, \theta )$, in the limit as $z \to \partial \Omega$). On the
other hand $f(z) \to 0$ and $\varphi (z)/f(z) \to 1$ as $z \to
\partial \Omega$. We may conclude that $k_g (\sigma ) \to -
4/(n+1)$ as $z \to \partial \Omega$. To complete the proof of
Klembeck's result we must compute the sectional curvature of the
$2$-plane $\sigma_0 \subset T_z (\Omega )$ spanned by $\{ N_z ,
T_z \}$ (remember that $J N = T$). Note first that
\[ N(f) = f^2 \left( \frac{2}{\varphi^2} + N(r) \right) . \]
Let us set for simplicity
\[ g = N(r) + \frac{4}{\varphi^2} - \frac{2r}{\varphi} \, , \;\;\;
h = N(r) + \frac{4}{\varphi^2} - \frac{6r}{\varphi} + 4 r^2 \, .
\]
We these notations let us recall that (by (\ref{e:b32}))
\begin{equation}
\nabla^g_T T = - \frac{1}{2} \; X_r - \frac{f}{2} \left\{ T(r) T -
h N \right\}
\end{equation} where $X_r = \nabla^H r$. Using also (\ref{e:b29}) for
$X = X_r$ we obtain
\[ - 2 \nabla^g_N \nabla^g_T T = \nabla_N X_r - \frac{1}{\varphi} \,
X_r + \frac{f}{2} \left\{ (\phi X_r )(r) T - X_r (r) N \right\} +
\] \[ +  N(f) \{ T(r) T - h N \} + f \left\{ N(T(r)) T + T(r) \nabla^g_N T - N(h) N - h \nabla^g_N
N \right\} .  \] Let us recall that (by (\ref{e:b30}) and
(\ref{e:b33}))
\begin{equation}
\nabla^g_N T = - \frac{1}{2} \, \phi \, X_r - \frac{f}{2} \left\{
g T + T(r) N \right\} ,
\end{equation}
\begin{equation}
\nabla^g_N N = - \frac{1}{2} \; X_r + \frac{f}{2} \left\{ T(r) T -
g  N \right\} .
\end{equation}
Using these identities and the expression of $N(f)$ gives (after
some simplifications)
\begin{equation}
- 2 \nabla^g_N \nabla^g_T T = \nabla_N X_r + \left( \frac{f h}{2}
- \frac{1}{\varphi} \right) X_r - \frac{f}{2} \, T(r) \, \phi X_r
+ \label{e:430}
\end{equation}
\[ + \frac{f}{2} \left\{ 2 f \left( \frac{2}{\varphi^2} + N(r)
\right) T(r) + 2 N(T(r)) - f(g+h) T(r) \right\} T - \]
\[ - \frac{f}{2} \left\{ g_\theta (X_r , X_r ) + 2 f h \left(
\frac{2}{\varphi^2} + N(r) \right) + 2 N(h) + f [T(r)^2 - gh]
\right\} N \] because of
\[ (\phi X_r )(r) = g_\theta (\nabla r , \phi X_r ) = g_\theta
(X_r , \phi X_r ) = 0, \]
\[ X_r (r) = g_\theta (\nabla^H r , X_r ) = g_\theta (X_r , X_r ).
\]
Similarly
\begin{equation}
-2 \nabla^g_T \nabla^g_N T = \nabla_T \phi X_r + \left(
\frac{1}{f} - \frac{fg}{2} \right) X_r + \frac{f}{2} \, T(r) \,
\phi X_r + \label{e:431}
\end{equation}
\[ + \frac{f}{2} \left\{ 2 T(g) + f (g-h) T(r)
\right\} T + \]
\[ + \frac{f}{2} \left\{ g_\theta (X_r , X_r ) + 2
T^2 (r) + f[T(r)^2 + gh] \right\} N . \] Here $T^2 (r) = T(T(r))$.
Let us set $\tau (W_\alpha ) = A_\alpha^{\overline{\beta}}
W_{\overline{\beta}}$. To compute the last term in the right hand
member of
\begin{equation}
R^g (N,T)T = \nabla^g_N \nabla^g_T T - \nabla^g_T \nabla^g_N T -
\nabla^g_{[N,T]} T \label{e:432}
\end{equation}
note first that $T(f) = f^2 \, T(r)$. On the other hand we may use
the decomposition (\ref{e:A.5}) so that
\[ \nabla^g_{[N,T]} T = r X_r + f r T(r) T - \frac{f}{2} \{
g_\theta (X_r , X_r ) + 2 r h \} N + \]
\[ + \left( i r^{\overline{\alpha}} A_{\overline{\alpha}}^\beta -
\frac{1}{f} r^\beta \right) W_\beta - \left( i r^\alpha
A_\alpha^{\overline{\beta}} + \frac{1}{f} r^{\overline{\beta}}
\right) W_{\overline{\beta}} \] (where
$A_{\overline{\alpha}}^\beta =
\overline{A_{\alpha}^{\overline{\beta}}}$) and by taking into
account that
\[ \left( i r^{\overline{\alpha}} A_{\overline{\alpha}}^\beta -
\frac{1}{f} r^\beta \right) W_\beta - \left( i r^\alpha
A_\alpha^{\overline{\beta}} + \frac{1}{f} r^{\overline{\beta}}
\right) W_{\overline{\beta}}  = - \frac{1}{f} \, X_r - \tau (\phi
X_r ) \] we may conclude that
\begin{equation}
\nabla^g_{[N,T]} T = \left( r - \frac{1}{f} \right) X_r - \tau
(\phi X_r ) + \label{e:433}
\end{equation}
\[ + f r T(r) T - \frac{f}{2} \{ g_\theta (X_r , X_r ) + 2 r h \}
N. \] Finally (by plugging into (\ref{e:432})  from
(\ref{e:430})-(\ref{e:431}) and (\ref{e:433}))
\begin{equation}
- 2 R^g (N , T) T = \nabla_N X_r - \nabla_T \phi X_r - f T(r) \phi
X_r - 2 \tau (\phi X_r ) + \label{e:434}
\end{equation}
\[ + \left( 2 r + \frac{f}{2} (g+h) - \frac{1}{\varphi} -
\frac{3}{f} \right) X_r + \]
\[ + f \left\{ f \left( \frac{2}{\varphi^2} + N(r) \right) T(r) +
N(T(r)) - T(g) + (2 r - fg) T(r) \right\} T - \]
\[ - f \left\{ 2 \| X_r \|^2 + f h \left(
\frac{2}{\varphi^2} + N(r) \right) + N(h) + f T(r)^2 + T^2 (r) + 2
r h \right\} N. \] Here $\| X_r \|^2 = g_\theta (X_r , X_r )$. Let
us take the inner product of (\ref{e:434}) with $N$ and use
(\ref{e:b5})-(\ref{e:b7}). We obtain
\[ 2 g(R^g (N,T)T,N) = \] \[ = \frac{n+1}{\varphi} \left\{ 2 \| X_r \|^2 +
f h \left( \frac{2}{\varphi^2} + N(r) \right) + \right. \] \[ +
\left.  N(h) + f T(r)^2 + T^2 (r) + 2 r h \right\} \] and dividing
by
\[ g(N,N) g(T,T) - g(N,T)^2 = \frac{1}{f^2} \left( \frac{n+1}{\varphi}
\right)^2 \] leads to
\[ 2 \frac{g(R^g (N,T)T , N)}{g(N,N) g(T,T) - g(N,T)^2} = \]
\[ = \frac{f^2 \varphi}{n+1} \left\{ 2 \| X_r \|^2 + T^2 (r) + f T(r)^2
+ 2 h r + N(h) + f h N(r) + 2 \frac{f h}{\varphi^2} \right\} .
\]
It remains that we perform an elementary asymptotic analysis of
the right hand member of the previous identity when $z \to
\partial \Omega$ (equivalently when $\varphi \to 0)$. As $r \in
C^\infty (\overline{\Omega})$ (cf. \cite{kn:LeMe}) the terms $\|
X_r \|^2$, $T^2 (r)$, $T(r)^2$ and $N(r)$ stay finite at the
boundary. Also (by recalling the expression of $h$) $f^2 \varphi h
\to 0$ as $\varphi \to 0$. Moreover
\[ 2 \frac{f^2 \varphi}{n+1} \frac{fh}{\varphi^2} = \frac{2}{n+1}
\frac{f}{\varphi} \left[ f^2 N(r) + \frac{4}{(1- r \varphi )^2} -
\frac{6 f^2 r}{\varphi} + 4 f^2 r^2 \right] \to \frac{8}{n+1} \, ,
\]
\[ N(h) = N^2 (r) + 4 N(r^2 ) - \frac{16}{\varphi^3} +
\frac{12r}{\varphi^2} - \frac{6}{\varphi} \, N(r) , \]
\[ \frac{f^2 \varphi}{n+1} \, N(h) \to - \frac{16}{n+1} \, , \]
as $\varphi \to 0$ hence
\[ k_g (\sigma_0 ) \to - \frac{4}{n+1} \, , \;\;\; z \to \partial
\Omega . \] Klembeck's theorem is proved.

\end{document}